\documentclass{amsart}
\usepackage{amssymb}
\usepackage{amscd}

\theoremstyle{plain}
\newtheorem{thm}{Theorem}[section]

\newtheorem{pro}[thm]{Proposition}

\theoremstyle{definition}

\theoremstyle{remark}

\newcommand{\C}{\mathbb{C}}

\newcommand{\Z}{\mathbb{Z}}

\newcommand{\PP}{\mathbb{P}}

\DeclareMathOperator{\Der}{Der}
\DeclareMathOperator{\Coh}{Coh}

\DeclareMathOperator{\Hom}{Hom}

\DeclareMathOperator{\End}{End}
\DeclareMathOperator{\Ext}{Ext}

\DeclareMathOperator{\Mod}{Mod}

\DeclareMathOperator{\Res}{Res}

\newcommand{\OO}{\mathcal{O}}
\newcommand{\F}{\mathbf{F}}
\newcommand{\E}{\mathbf{E}}

\newcommand{\XX}{\mathbb{X}}

\newcommand{\ox}{\vec{x}}

\newcommand{\oc}{\vec{c}}
\newcommand{\oom}{\vec{\omega}}

\begin{document}
\title{Connections for weighted projective lines}

\author{William Crawley-Boevey}
\address{Department of Pure Mathematics, University of Leeds, Leeds LS2 9JT, UK}
\email{w.crawley-boevey@leeds.ac.uk}


%14  (1940-now) Algebraic geometry
%14H  (1973-now) Curves
%14H45  (1973-now) Special curves and curves of low genus

%16  (1959-now) Associative rings and algebras [For the commutative case, see 13-XX]
%16G  (1991-now) Representation theory of rings and algebras
%16G20  (1991-now) Representations of quivers and partially ordered sets

\begin{abstract}
We introduce a notion of a connection on a coherent sheaf on a
weighted projective line (in the sense of Geigle and Lenzing).
Using a theorem of H\"ubner and Lenzing we show, under a mild hypothesis,
that if one considers coherent sheaves equipped with such a connection,
and one passes to the perpendicular category to a nonzero vector bundle without self-extensions,
then the resulting category is equivalent to the category of representations of a deformed preprojective algebra.
\end{abstract}
\maketitle

\section{Introduction}
Consider the weighted projective line $\XX$ consisting of the complex projective line $X = \PP^1$,
a collection $D = (a_1,\dots,a_k)$ of
distinct points of $X$, and a weighting $\mathbf{w} = (w_1,\dots,w_k)$, with the $w_i$
being integers $\ge 2$ (or also 1, which is equivalent to the point not being marked).

Weighted projective lines, and the category $\Coh\XX$ of coherent sheaves over $\XX$,
were introduced by Geigle and Lenzing \cite{GL}.
Later, Lenzing \cite{L} showed that the category of vector bundles $\E$ on $\XX$
is equivalent to the category of vector
bundles $E$ on $\PP^1$ equipped with a (quasi)parabolic structure of
weight type $(D,\mathbf{w})$, that is, flags of subspaces
\[
E_{a_i} = E_{i0} \supseteq E_{i1} \supseteq \dots \supseteq E_{i,w_i-1} \supseteq E_{i,w_i}=0
\]
of the fibres of $E$ at each marked point.

It is natural to consider connections $\nabla:E\to E\otimes \Omega^1_X(\log D)$
on a vector bundle on $\PP^1$, analytic except possibly for logarithmic poles on $D$,
and their residues
\[
\Res_{a_i} \nabla:E_{a_i}\to E_{a_i}.
\]
Given a collection of complex numbers $\zeta = (\zeta_{is})$ ($1\le i\le k$, $1\le s\le w_i$),
if $E$ is a parabolic bundle on $\PP^1$ of weight type $(D,\mathbf{w})$,
then $\nabla$ is said to be a $\zeta$-connection \cite[\S7]{CBipb} provided that
\[
(\Res_{a_i} \nabla - \zeta_{is} 1)(E_{i,s-1}) \subseteq E_{is}
\]
for all $1\le i\le k$ and $1\le s\le w_i$.
This condition has the effect of fixing the eigenvalues of the residues, and if one knows the dimensions
of the subspaces in the flags, it also fixes the Jordan block structure.

In this first part of this article we define the appropriate notion of a $\zeta$-connection for any
coherent sheaf on $\XX$, not just those corresponding to parabolic bundles.

H\"ubner and Lenzing \cite{HL} proved that if $\F$ is a nonzero vector bundle on $\XX$ 
without self-extensions, then the perpendicular category
\[
\F^\perp = \{ \E\in\Coh\XX \mid \Hom(\F,\E) = \Ext^1(\F,\E) = 0 \}
\]
is equivalent to the category of representations of a finite-dimensional hereditary algebra $A$.
In the second part of this article we show that, under a mild hypothesis, the category
of coherent sheaves in $\F^\perp$ equipped with
a $\zeta$-connection is equivalent to the category of representations of a
deformed preprojective algebra $\Pi^\lambda(A)$ in the sense of \cite{CBconze}.
These results will be used elsewhere in our work on the Deligne-Simpson problem,
see \cite{CBicm} for an overview.

I would like to thank H. Lenzing for a number of very helpful discussions.

\section{Weighted projective lines via patching}
Lenzing \cite[\S4.2]{L} showed how one can define the category $\Coh\XX$ for one
marked point in terms of $n$-cycles, and then an iterative construction generalizes this to any number of marked points.
Here we outline a variation which seems more direct. Everything is implicit in \cite{GL,L}, so we omit details.

Let $A$ be a commutative integral domain and $\mathfrak{m}$ an invertible ideal in $A$.
If $M$ is an $A$-module, we define $M(s) = \mathfrak{m}^{-s}\otimes M$ for $s\in\Z$.
Since $\mathfrak{m}$ is invertible, we can identify $M(s)(t)$ with $M(s+t)$ for any $s,t$.
Moreover for $s\le t$ the inclusion $\mathfrak{m}^{-s}\subseteq \mathfrak{m}^{-t}$  induces a natural map $M(s)\to M(t)$.

If $n\ge 1$, by an \emph{$n$-cycle of $A$-modules concentrated at $\mathfrak{m}$}
one means a collection $E = (E_s,\phi_s)_{s\in\Z}$ of $A$-modules and homomorphisms
\[
\dots \to E_2 \xrightarrow{\phi_{1}} E_1 \xrightarrow{\phi_{0}} E_0 \xrightarrow{\phi_{-1}} E_{-1} \to \dots
\]
such that for all $s\in\Z$ we have $E_{s-n} = E_s(1)$, and the composition
\[
E_{s} \xrightarrow{\phi_{s-1}} E_{s-1} \xrightarrow{\phi_{s-2}} \dots \xrightarrow{\phi_{s-n}} E_{s-n} = E_s(1)
\]
is the natural map. There is a natural category of such cycles, in which the morphisms $u:E\to F$ are
given by $A$-module maps $u_s:E_s\to F_s$ which are \emph{$n$-periodic}, that is, for all $s\in\Z$, the map
$u_{s-n}:E_{s-n}\to F_{s-n}$ is the same as the map $E_s(1)\to F_s(1)$ induced by $u_s$, and
satisfy $\phi_s u_{s+1} = u_s \phi_s$ for all $s$.
This category is easily seen to be equivalent to the category of modules
for the ring $T_n(A,\mathfrak{m})$ of $n\times n$ matrices of shape
\[
T_n(A,\mathfrak{m}) = \begin{pmatrix}
A & A & A & \hdots \\
\mathfrak{m} & A & A & \\
\mathfrak{m} & \mathfrak{m} & A & \\
\vdots & & & \ddots
\end{pmatrix}.
\]

Now fix marked points $D = (a_1,\dots,a_k)$ and a weighting $\mathbf{w} = (w_1,\dots,w_k)$ as above.
By adding additional marked points, with trivial weighting~1, we may assume that $k\ge 2$.
Thus we can fix an affine open covering $X_i$ ($1\le i \le k$) of $X=\PP^1$,
such that each $X_i$ contains $a_i$, but no other marked points. Let $A_i$ be the
coordinate ring of $X_i$, and let $\mathfrak{m}_i$
be the maximal ideal corresponding to the point~$a_i$.

Recall that $\Coh X$ can be identified with the category whose objects are collections $(E^i_0,g_{ij})$
where each $E^i_0$ ($1\le i \le k$) is a finitely generated $A_i$-module,
and the $g_{ij}$ ($1\le i,j\le k$) are patching isomorphisms
\[
g_{ij}:A_{ij} \otimes_{A_j} E^j_0 \to A_{ij} \otimes_{A_i} E^i_0,
\]
where $A_{ij}$ is the coordinate ring of $X_i\cap X_j$, subject to the compatibility conditions
\[
g_{ij} g_{j\ell } = g_{i\ell}
\]
for all $i,j,\ell$, for modules over the coordinate ring of $X_i\cap X_j \cap X_k$.
The morphisms from $(E^i_0,g_{ij})$ to $(F^i_0,h_{ij})$ in this category
are collections $f = (f^i_0)_{1\le i\le k}$ where each $f^i_0:E^i_0\to F^i_0$
is a morphism of $A_i$-modules, and compatible in the sense the square
\[
\begin{CD}
A_{ij} \otimes_{A_j} E^j_0 @>1\otimes f^j_0>> A_{ij} \otimes_{A_j} F^j_0 \\
@Vg_{ij}VV @Vh_{ij}VV \\
A_{ij} \otimes_{A_i} E^i_0 @>1\otimes f^i_0>> A_{ij} \otimes_{A_i} F^i_0
\end{CD}
\]
commutes for all $i,j$.

The category $\Coh\XX$ is defined as follows.
Objects are collections $\E = (E^i,g_{ij})$
where each $E^i = (E^i_s,\phi^i_s)$ ($1\le i \le k$)
is a $w_i$-cycle of finitely generated $A_i$-modules concentrated at $\mathfrak{m}_i$,
and the $g_{ij}$ ($1\le i,j\le k$) are patching isomorphisms
\[
g_{ij}:A_{ij} \otimes_{A_j} E^j_0 \to A_{ij} \otimes_{A_i} E^i_0,
\]
where $A_{ij}$ is the coordinate ring of $X_i\cap X_j$, subject to the compatibility conditions
\[
g_{ij} g_{j\ell } = g_{i\ell}
\]
for all $i,j,\ell$, for modules over the coordinate ring of $X_i\cap X_j \cap X_k$.
Morphisms $\E\to \F$, where $\E=(E^i,g_{ij})$ and $\F=(F^i,h_{ij})$,
are collections $f = (f^i)_{1\le i\le k}$ where each $f^i:E^i\to F^i$
is a morphism of $w_i$-cycles of $A_i$-modules concentrated at $\mathfrak{m}_i$, so $f^i = (f^i_s)_{s\in\Z}$,
and they are compatible in the sense the square
\[
\begin{CD}
A_{ij} \otimes_{A_j} E^j_0 @>1\otimes f^j_0>> A_{ij} \otimes_{A_j} F^j_0 \\
@Vg_{ij}VV @Vh_{ij}VV \\
A_{ij} \otimes_{A_i} E^i_0 @>1\otimes f^i_0>> A_{ij} \otimes_{A_i} F^i_0
\end{CD}
\]
commutes for all $i,j$.

Clearly there is a forgetful functor $\Coh\XX\to\Coh X$ which replaces each cycle $E^i$ by the module $E^i_0$.

By definition a \emph{vector bundle} on $\XX$ is a coherent sheaf $\E = (E^i,g_{ij})$ with the property
that for all $i$ and $s$, the modules $E^i_s$ occurring in the cycles $E^i = (E^i_s,\phi^i_s)$
are projective $A_i$-modules.
If $\E = (E^i,g_{ij})$ is a vector bundle and $E$ is the corresponding vector bundle on $X$,
then the fibre of $E$ at $a_i$ is
\[
E_{a_i} = E^i_0 \otimes_{A_i} A_i/\mathfrak{m}_i
\]
and there is a flag of subspaces
\[
E_{a_i} = E_{i0} \supseteq E_{i1} \supseteq \dots \supseteq E_{i,w_i-1} \supseteq E_{i,w_i}=0,
\]
where $E_{is}$ is the image of the composition of $\phi^i_0\phi^1_1 \cdots\phi^i_{s-1} : E^i_s \to E^i_0$ with the projection $E^i_0\to E_{a_i}$.
This construction gives the equivalence between the the category of vector bundles on $\XX$ and the category of
parabolic bundles on $X$ of weight type $(D,\mathbf{w})$.

Given integers $r_i$ for $1\le i \le k$, the shift $\E(\sum_{i=1}^k r_i\vec x_\ell)$
of a coherent sheaf $\E = (E^i,g_{ij})$ with $E^i = (E^i_s,\phi^i_s)$
is defined to be $\F = (F^i,h_{ij})$, where $F^i = (F^i_s,\chi^i_s)$ with $F^i_s = E^\ell_{s-r_i}$
and $\chi^i_s = \phi^i_{s-r_i}$,
and the $h_{ij}$ are obtained from the $g_{ij}$ by conjugating, where necessary, using the maps
\[
\dots \to
A_{ij}\otimes_{A_i} E^i_2 \xrightarrow{1\otimes \phi^i_1}
A_{ij}\otimes_{A_i} E^i_1 \xrightarrow{1\otimes \phi^i_0}
A_{ij}\otimes_{A_i} E^i_0 \xrightarrow{1\otimes \phi^i_{-1}}
A_{ij}\otimes_{A_i} E^i_{-1} \to \dots\ .
\]
These are all isomorphisms for $i\neq j$, since the condition that $a_i\notin X_j$ implies that $A_{ij}\mathfrak{m}_i = A_{ij}$.
For example, if all $r_i>0$ then $h_{ij}$ ($i\neq j$) is the composition
\[
A_{ij}\otimes F^j_0 = A_{ij}\otimes E^j_{-r_j}
\xrightarrow{p^{-1}}
A_{ij}\otimes E^j_0
\xrightarrow{g_{ij}}
A_{ij}\otimes E^i_0
\xrightarrow{q}
A_{ij}\otimes E^i_{-r_i} = A_{ij}\otimes F^i_0,
\]
where $p = (1\otimes \phi^j_{-r_j})\cdots(1\otimes \phi^j_{-1})$
and $q = (1\otimes \phi^i_{-r_i})\cdots(1\otimes \phi^i_{-1})$.

If $\E$ is a coherent sheaf on $\XX$ and $L$ is a coherent sheaf on $X$, then
then there is a coherent sheaf $\E\otimes L$ on $\XX$. If $\E = (E^i,g_{ij})$, $E^i = (E^i_s,\phi^i_s)$, and
the restriction of $L$ to $X_i$ is given by an $A_i$-module $L_i$, then
$\E\otimes L$ is given by the modules $E^i_s\otimes_{A_i} L_i$.
In this terminology, the expression $(\E\otimes \Omega^1_X)(\sum_{i=1}^k (w_i-1)\vec x_i)$
can be identified with the twist $\E(\vec \omega)$, as in \cite{GL}.

\section{Connections for weighted projective lines}
Atiyah \cite[\S4]{Atiyah} observed that a connection on a vector bundle $E$ on a
complex manifold $X$ is the same thing as a section of a certain exact sequence
\[
\mathcal{B}(E) : 0\to E\otimes\Omega^1_X \to \mathrm{D}(E) \to E\to 0
\]
Similarly, Mihai \cite{Mihai1,Mihai2} identified logarithmic connections on $E$ with sections of the pushout sequence
\[
\mathcal{B}_D(E) : 0\to E\otimes\Omega^1_X(\log D) \to \mathrm{D}_D(E) \to E\to 0.
\]
Here we put this in the setting of weighted projective lines.

\begin{thm}
Let $\mathbb{X}$ be a weighted projective line over $\C$ associated to $D$ and $\mathbf{w}$.
Fix a collection of complex numbers $\zeta = (\zeta_{is})$ ($1\le i\le k$, $1\le s\le w_i$).
Then there are exact sequences
\[
\mathcal{B}_\zeta(\E) : 0\to \E(\oom) \to \mathbf{D}_\zeta(\E) \to \E \to 0
\]
defined for any coherent sheaf $\E$ on $\XX$, and functorial in $\E$,
whose sections, in case $\E$ is a vector bundle,
are in 1-1 correspondence with the $\zeta$-connections on the corresponding parabolic bundle.
\end{thm}

In view of this, we define a \emph{$\zeta$-connection} on $\E$ to be a section of $\mathcal{B}_\zeta(\E)$.

\begin{proof}
Since $X=\PP^1$, an easy calculation shows that given any $c_1,\dots,c_k\in\C$
with $\sum_i c_i=0$, there is a global section of $\Omega^1_X(\log D)$ whose
residue at $a_i$ is $c_i$.
For any $E$ this gives a homomorphism $E\to E\otimes \Omega^1_X(\log D)$, and we
can add this homomorphism to a logarithmic connection $\nabla:E\to E\otimes \Omega^1_X(\log D)$,
changing a $\zeta$-connection into a $\zeta'$-connection, where $\zeta'_{is}=\zeta_{is}+c_i$.
By making such a change we may assume that $\zeta_{i1}$ is independent of $i$. Denote the common value by $\zeta_0$.

Consider an affine open piece $X_i$.
We identify coherent sheaves on $X_i$ with the corresponding finitely generated $A_i$-module.

We write the restriction of $\Omega^1_X(\log D)$ to $X_i$ as $\Omega^1_{X_i}(\log D)$,
but as $a_i$ is the only point of $D$ in $X_i$, it is the same as $\Omega^1_{X_i}(\log a_i)$,
that is, $\mathfrak{m}_i^{-1}\otimes \Omega^1_{X_i}$. Let
\[
\rho_{ij}:A_{ij}\otimes_{A_j} \Omega^1_{X_j}(D) \to A_{ij}\otimes_{A_i} \Omega^1_{X_i}(D)
\]
be the patching isomorphisms for $\Omega^1_X(\log D)$.

Recall that $\E$ is given on $X_i$ by a $w_i$-cycle $E^i = (E^i_s,\phi^i_s)_{s\in\Z}$.
Then $\E(\oom)$ is given by the cycle whose $s$th term is
\[
E^i_{s-(w_i-1)}\otimes \Omega^1_{X_i}
\cong E^i_{s+1}\otimes \mathfrak{m}_i^{-1}\otimes \Omega^1_{X_i}
\cong E^i_{s+1}\otimes \Omega^1_{X_i}(\log D)
\]

For $s\in\Z$, we define
\[
D^i_s = E^i_s \ltimes (E^i_{s+1}\otimes \Omega^1_{X_i}(\log D))
\]
where the notation $\ltimes$ means the direct sum as vector spaces,
considered as an $A_i$-module via the action with
\[
a(e,f\otimes \gamma) = (ae,\phi^i_{s+1}\phi^i_{s+2}\dots\phi^i_{s+w_i-1}(x_i\otimes e)\otimes \frac{da}{x_i} + af\otimes \gamma)
\]
for $a\in A_i$, $e\in E^i_s$, $f\in E^i_{s+1}$ and $\gamma\in \Omega^1_{X_i}(\log D)$.
Here $x_i$ is a local coordinate at $a_i$, so a generator of $\mathfrak{m}_i$ as an ideal in $A_i$,
$x_i\otimes e$ is in $\mathfrak{m}_i\otimes E^i_s = E^i_s(-1) = E^i_{s+w_i}$, and the composition $\phi_{s+1}\phi_{s+2}\dots\phi_{s+w_i-1}$
sends it to an element of $E^i_{s+1}$.

Then the canonical inclusion and projection define exact sequences of $A_i$-modules
\[
\mathcal{B}^i_s(\E) : 0 \to E^i_{s+1}\otimes \Omega^1_{X_i}(\log D) \to D^i_s \to E^i_s \to 0.
\]
We define maps $\psi^i_s:D^i_{s+1}\to D^i_s$, that is,
\[
\psi^i_s : E^i_{s+1} \ltimes (E^i_{s+2}\otimes \Omega^1_{X_i}(\log D)) \to E^i_s \ltimes (E^i_{s+1}\otimes \Omega^1_{X_i}(\log D))
\]
by
\[
\psi^i_s = \begin{pmatrix} \phi^i_s & 0 \\ 1\otimes (\zeta_{i,s+2(w_i)} - \zeta_{i,s+1(w_i)}) \frac{dx_i}{x_i} & \phi^i_{s+1}\otimes 1 \end{pmatrix}.
\]
where if $j\in \Z$ then $j(w_i)$ denotes the integer in the range $1,\dots,w_i$ congruent to $j$ modulo~$w_i$.
It is easy to see that they turn the $D^i_s$ ($s\in\Z$) into a $w_i$-cycle of $A_i$-modules
concentrated at $\mathfrak{m}_i$, and that there are morphisms of exact sequences
\[
\begin{CD}
\mathcal{B}^i_{s+1}(\E) : 0
@>>> E^i_{s+2} \otimes \Omega^1_{X_i}(\log D)
@>>> D^i_{s+1}
@>>> E^i_{s+1}
@>>> 0
\\
& & @V \phi^i_{s+1}\otimes 1 VV @V\psi^i_s VV @V \phi^i_s VV
\\
\mathcal{B}^i_s(\E) : 0
@>>> E^i_{s+1} \otimes \Omega^1_{X_i}(\log D)
@>>> D^i_s
@>>> E^i_s
@>>> 0.
\end{CD}
\]

We next define patching isomorphisms to turn the $D^i_s$ into a
coherent sheaf $\mathbf{D}_\zeta(\E)$ on $\mathbb{X}$.
Let $g_{ij}:A_{ij} \otimes_{A_j} E^j_0 \to A_{ij} \otimes_{A_i} E^i_0$,
be the patching isomorphisms for $\E$.
The image $E$ of $\E$ under the forgetful functor is defined by the $A_i$-modules $E^i_0$
and the patching isomorphisms $g_{ij}$.
Moreover $E\otimes \Omega^1_{X} (\log D)$ is given by the modules $E^i_0 \otimes_{A_i} \Omega^1_{X_i} (\log D)$
with patching isomorphisms $g_{ij}\otimes \rho_{ij}$, which can be considered as mappings
\[
A_{ij} \otimes_{A_j} E^j_0\otimes_{A_j} \Omega^1(\log D)
\to
A_{ij} \otimes_{A_i} E^i_0\otimes_{A_j} \Omega^1(\log D).
\]
The sheaf $\mathrm{D}_D(E)$ in the Mihai sequence is obtained from the $A_i$-modules
\[
\mathrm{D}_D(E^i_0) = E^i_0 \ltimes (E^i_0 \otimes \Omega^1_{X_i}(\log D))
\]
with action
\[
a (e,f\otimes \gamma) = (ae, e\otimes da+af\otimes \gamma)
\]
and suitable patching isomorphisms
\[
\ell_{ij} : A_{ij} \otimes_{A_j} \mathrm{D}_D(E^j_0) \to A_{ij} \otimes_{A_i} \mathrm{D}_D(E^i_0).
\]
For $i\neq j$ the mapping
\[
1\otimes \phi^j_0 \otimes 1:
A_{ij} \otimes_{A_j} E^j_1\otimes_{A_j} \Omega^1(\log D)
\to
A_{ij} \otimes_{A_j} E^j_0\otimes_{A_j} \Omega^1(\log D),
\]
is invertible. The patching isomorphisms for $\E(\oom)$ are
\[
h_{ij} =
\begin{cases}
(1\otimes \phi^i_0 \otimes 1)^{-1} \circ (g_{ij}\otimes \rho_{ij}) \circ (1\otimes \phi^j_0 \otimes 1)
& (i\neq j)
\\
1 & (i=j).
\end{cases}
\]

There is a morphism of exact sequences of $A_i$-modules
\[
\begin{CD}
\mathcal{B}^i_0(\E) : 0
@>>> E^i_1\otimes \Omega^1(\log D)
@>>> D^i_0
@>>> E^i_0
@>>> 0
\\
& & @V\phi^i_0\otimes 1 VV @V\alpha_i VV @|
\\
\mathcal{B}(E^i_0) : 0
@>>> E^i_0\otimes \Omega^1(\log D)
@>>> \mathrm{D}_D(E^i_0)
@>>> E^i_0
@>>> 0
\end{CD}
\]
where
\[
\alpha_i : E^i_0 \ltimes (E^i_1\otimes \Omega^1_{X_i}(\log D))
\to
E^i_0 \ltimes (E^i_0\otimes \Omega^1_{X_i}(\log D)).
\]
is the mapping with matrix
\[
\begin{pmatrix} 1 & 0 \\ 1\otimes \zeta_{i1} \frac{dx_i}{x_i} & \phi^i_0\otimes 1 \end{pmatrix}.
\]
If $i\neq j$ then the mapping
\[
1\otimes \alpha_i : A_{ij}\otimes_{A_i} D^i_0 \to A_{ij}\otimes_{A_i} \mathrm{D}_D(E^i_0)
\]
is an isomorphism. The patching maps for $\mathbf{D}_\zeta(\E)$ are then the maps
\[
k_{ij} : A_{ij}\otimes_{A_j} D^j_0 \to A_{ij}\otimes_{A_i} D^i_0
\]
defined by
\[
k_{ij} =
\begin{cases}
(1\otimes \alpha_i)^{-1} \ell_{ij} (1\otimes \alpha_j)
& (i\neq j)
\\
1 & (i=j).
\end{cases}
\]
Clearly the exact sequences above define exact sequences
\[
0\to \E(\oom) \to \mathbf{D}_\zeta(\E) \to \E \to 0
\]
of coherent sheaves on $\XX$, and this is all functorial in $\E$.

Now suppose that $\E$ is a vector bundle and consider a section $\sigma$ of the sequence $\mathcal{B}_\zeta(\E)$.
This is given by sections $\sigma^i_s$ of $\mathcal{B}^i_s(\E)$
which are $w_i$-periodic in $s$,
satisfy
\[
\psi^i_s \sigma^i_{s+1} = \sigma^i_s \phi^i_s,
\]
and commutativity of the squares
\[
\begin{CD}
A_{ij}\otimes_{A_j} E^j_0 @>g_{ij}>> A_{ij}\otimes_{A_i} E^i_0 \\
@V 1\otimes \sigma^j_0 VV @V 1\otimes \sigma^i_0 VV \\
A_{ij}\otimes_{A_j} D^j_0 @>k_{ij}>> A_{ij}\otimes_{A_i} D^i_0.
\end{CD}
\]
Now one can write
\[
\sigma^i_s = \begin{pmatrix} 1 \\ \nabla^i_s \end{pmatrix} : E^i_0 \to D^i_s = E^i_s \ltimes (E^i_{s+1} \otimes \Omega^1_{X_i}(\log D))
\]
for some mappings $\nabla^i_s$, and then
\[
\nabla^i_s \phi^i_s = (\phi^i_{s+1}\otimes 1)\nabla^i_{s+1} + 1 \otimes (\zeta_{i,s+2(w_i)}-\zeta_{i,s+1(w_i)})\frac{dx_i}{x_i}
\]

The composition $\alpha_i \sigma^i_0$ is a section of the Mihai sequence for $E^i_0$, so of the form
\[
\alpha_i\sigma^i_0 = \begin{pmatrix} 1 \\ \nabla^i \end{pmatrix} : E^i_0 \to E^i_0 \ltimes (E^i_0 \otimes \Omega^1_{X_i}(\log D))
\]
for some connection $\nabla^i = \alpha_i \sigma^i_0 : E^i_0 \to E^i_0\otimes \Omega^1_{X_i}(\log D)$.
Thus
\[
\nabla^i = (\phi^i_0\otimes 1)\nabla^i_0 + 1\otimes \zeta_{i1} \frac{dx_i}{x_i}.
\]
Moreover the definition of $k_{ij}$ ensures that the $\alpha_i \sigma^i_0$ patch together
to give a section of $\mathcal{B}_D(E)$, so the $\nabla^i$ give a logarithmic connection
\[
\nabla : E \to E\otimes \Omega^1_{X}(\log D).
\]

The residue at $a_i$ defines an $A_i$-module map
\[
\Res_i : \Omega^1_{X_i}(\log D) = \mathfrak{m}_i^{-1} \otimes_{A_i} \Omega^1_{X_i} \to A_i/\mathfrak{m}_i.
\]
If $u\in \mathfrak{m}_i^{-1}$ and $m \in \mathfrak{m}_i$, this map sends $u \otimes dm \in \mathfrak{m}_i^{-1} \otimes_{A_i} \Omega^1_{X_i}$
to the image of $um$ in $A_i/\mathfrak{m}_i$.
The fibre of $E$ at $a_i$ is the vector space
\[
E_{a_i} = E^i_0 \otimes_{A_i} A_i/\mathfrak{m}_i.
\]
We denote by $\pi^i_s$ the projection $E^i_s \to E^i_s\otimes A_i/\mathfrak{m}_i$.
The composition
\[
(1\otimes \Res_i)\nabla^i : E^i_0 \to E^i_0 \otimes_{A_i} A_i/\mathfrak{m}_i.
\]
factors through $\pi^i_0$, and the residue of $\nabla$ at $a_i$ is the linear map $\Res_{a_i}\nabla \in \End(E_{a_i})$ with
\[
(\Res_{a_i}\nabla) \pi^i_0 = (1\otimes \Res_i)\nabla^i.
\]
Let
\[
R^i_s : E^i_s\otimes A_i/\mathfrak{m}_i \to E^i_{s+1}\otimes A_i/\mathfrak{m}_i
\]
be the map with $R^i_s \pi^i_s = (1\otimes \Res_i)\nabla^i_s$.
The formula relating $\nabla^i$ with $\nabla^i_0$ gives
\[
(\phi^i_0\otimes 1)R^i_0 = \Res_{a_i}\nabla  - \zeta_{i1}1
\]
which shows that $\Res_{a_i}\nabla - \zeta_{i1}1$ has image contained in $E_{i1}$,
and the formula relating $\nabla^i_s$ with $\nabla^i_{s+1}$ gives
\[
R^i_s (\phi^i_s\otimes 1) = (\phi^i_{s+1}\otimes 1)R^i_{s+1} + (\zeta_{i,s+2(w_i)}-\zeta_{i,s+1(w_i)})1.
\]
so that by induction, for $1\le s<w_i$,
\[
(\phi^i_0 \phi^1_i \dots \phi^i_s \otimes 1)R^i_s = (\Res_{a_i}\nabla - \zeta_{i,s+1}1)(\phi^i_0 \phi^1_1 \dots \phi^i_{s-1} \otimes 1).
\]
so that the restriction of $\Res_{a_i}\nabla - \zeta_{i,s+1}$ to $E_{is}$ has image contained in $E_{i,s+1}$.
This shows that $\nabla$ is a $\zeta$-connection on $E$.

Conversely, if $\nabla$ is a $\zeta$-connection on $E$, one can reverse this
construction and show that $\nabla$ arises from a section of $\mathcal{B}_\zeta(\E)$.
\end{proof}

\section{Perpendicular categories}
If $\F$ is a coherent sheaf on $\XX$, the perpendicular category $\F^\perp$ is the full subcategory
of $\Coh \XX$ with objects
\[
\F^\perp = \{ \E\in\Coh\XX \mid \Hom(\F,\E) = \Ext^1(\F,\E) = 0 \}.
\]
We begin with a theorem of H\"ubner and Lenzing \cite{HL}.
We include a copy of their proof, since the manuscript is unpublished.
We use the setup of \cite{LP}, and let $\Lambda$ be the canonical algebra associated to $\XX$.
We identify the indecomposable projective $\Lambda$-modules with the line bundles $\OO(\ox)$,
$0\le \ox\le \oc$, and we identify the category $\Mod_\ge (\Lambda)$ with the subcategory $\Coh_\ge \XX$ of $\Coh \XX$.

\begin{pro}[H\"ubner and Lenzing]
If $\F$ is a nonzero bundle on $\XX$ satisfying
\begin{enumerate}
\item$\Ext^1(\Lambda,\F)=0$, that is, $\F\in \Mod_+(\Lambda)$, and
\item$\mu \F > 2p + \delta(\oom)$, and
\end{enumerate}
then $\F^\perp \subseteq \Coh_{\ge}\XX$, and it coincides with
\[
\F^\perp_{\Mod(\Lambda)} = \{ M\in \Mod(\Lambda) \mid \Hom(\F,M) = \Ext^1(\F,M) = 0 \}.
\]
\end{pro}

\begin{proof}
First, if $M\in \Mod_-(\Lambda)$ then by \cite[Theorem 3.2]{LP}, $M = \E[1]$ for some
vector bundle $\E$ on $\XX$ with $\Hom(\Lambda,\E)=0$. Now
$\mu \E = \mu \E - \mu \OO \le p+\delta(\oom)$ by \cite[Theorem 2.7]{LP},
and $\mu \F(\oom) \ge 2p+2\delta(\oom)$ by (1), so
so $\mu \F(\oom) - \mu \E > p+\delta(\oom)$, and hence
$\Hom(\E,\F(\oom))\neq 0$ by \cite[Theorem 2.7]{LP}.
Thus by Serre duality $\Ext^1(\F,\E)\neq 0$, that is $\Hom(\F,M)\neq 0$.
Thus $M\notin \F^\perp$.
This shows that $\F^\perp_{\Mod(\Lambda)} \subseteq \Mod_{\ge}(\Lambda)$ and
consequently $\F^\perp_{\Mod(\Lambda)} \subseteq \F^\perp$.

Second, suppose that $\E$ is a vector bundle which is not in $\Mod_+(\Lambda)$.
Then $\Ext^1(\OO(\oom),\E)\neq 0$, so $\Hom(\E,\OO(\oc+\oom))\neq 0$.
Now $\mu \OO(\oc)=p$, so $\mu\F-\mu\OO(\oc) > p+\delta(\oom)$ by (2),
so $\Hom(\OO(\oc),\F)\neq 0$ by \cite[Theorem 2.7]{LP}.
Thus $\Hom(\OO(\oc+\oom),\F(\oom))\neq 0$. Since $\F$ is a vector bundle,
and $\OO(\oc+\oom)$ is a line bundle there must be a monomorphism
$\OO(\oc+\oom)\to \F(\oom)$. Its composition with a nonzero morphism
$\E\to \OO(\oc+\oom)$ is again nonzero, so $\Hom(\E,\F(\oom))\neq 0$,
and hence $\Ext^1(\F,\E)\neq 0$. Thus $\E\notin \F^\perp$.
This shows that $\F^\perp \subseteq \Mod_{\ge}(\Lambda)$ and
consequently $\F^\perp \subseteq \F^\perp_{\Mod(\Lambda)}$.
\end{proof}

\begin{thm}[H\"ubner and Lenzing]
\label{t:hltheorem}
If $\F$ is a nonzero vector bundle on $\XX$ with $\Ext^1(\F,\F)=0$, then
$\F^\perp$ is equivalent to the category of representations of a finite-dimensional hereditary algebra $A$,
and the inclusion functors $j_XX:\F^\perp\to \Coh\XX$ and $j_\Lambda:\F^\perp \to \Mod(\Lambda)$ are exact and
admit left adjoints $\ell_\XX$ and $\ell_\Lambda$. Moreover $\F$ can be extended to a tilting bundle
$\mathbf{T} = \F\oplus \mathbf{C}$ on $\XX$ such that $\mathbf{C}$ is projective in~$\F^\perp$.
\end{thm}

\begin{proof}
If $\ox\in \mathbf{L}(\mathbf{p})$, the twist $\E \mapsto \E(\ox)$
is an equivalence of categories, so it defines an equivalence $\F^\perp \mapsto \F(\ox)^\perp$.
Thus, we may replace $\F$ by any twist, and hence we may assume conditions (1) and (2) in the last proposition hold,
and also (3) $\Hom(\F,\Lambda) = 0$.

Since $\F\in\Mod_+(\Lambda)$, it has projective dimension at most 1 as a $\Lambda$-module by \cite[Theorem 3.2]{LP}.
The fact that the inclusion functors are exact and have left adjoints then follows from \cite[Proposition 3.2]{GLperp}.
It follows that $\ell_\Lambda(\Lambda)$ is a projective generator for $\F^\perp$, so it is equivalent to the category
of modules for some finite-dimensional algebra $A$. Moreover, as an abelian exact subcategory of the
hereditary category $\Coh\XX$, the category $\F^\perp$, and hence $A$, is also hereditary.

Finally, because of condition (3) and the construction in the proof
of \cite[Proposition 3.2]{GLperp}, $\ell_\Lambda(\OO(\ox))$ is the middle term of
the universal exact sequence
\[
0\to \OO(\ox) \to \ell_\Lambda(\OO(\ox)) \to \F\otimes\Ext^1(\F,\OO(\ox))\to 0.
\]
Thus $\ell_\Lambda(\OO(\ox))$ is a vector bundle and $\Ext^1(\ell_\Lambda(\OO(\ox)),\F)=0$.
One can thus take
\[
\mathbf{C} = \ell_\Lambda(\Lambda) = \bigoplus_{0\le \ox\le c} \ell_\Lambda(\OO(\ox)).
\]
\end{proof}

\section{Deformed preprojective algebras}
We now study perpendicular categories for coherent sheaves on a weighted projective line $\XX$
equipped with a $\zeta$-connection.

The deformed preprojective algebra $\Pi^\lambda(A)$ associated to a hereditary algebra $A$ is defined in \cite{CBconze}.
In the theorem of H\"ubner and Lenzing we may take $A$ to be basic, in which case it is isomorphic to
the path algebra of a quiver $Q$ without oriented cycles, and then
$\Pi^\lambda(A)$ is isomorphic to the deformed preprojective algebra $\Pi^\lambda(Q)$
as introduced by Crawley-Boevey and Holland \cite{CBH}.

\begin{thm}
Let $\F$ be a nonzero vector bundle on $\XX$ with $\Ext^1(\F,\F)=0$ and fix $\zeta$.
Let $\F^\perp$ be equivalent to $\Mod(A)$ for some finite-dimensional hereditary algebra $A$,
and assume that $A$ has no non-zero projective-injective module.
Then the category of coherent sheaves in $\F^\perp$ equipped with a
$\zeta$-connection is equivalent to
the category of representations of the deformed preprojective algebra $\Pi^\lambda(A)$
for some \mbox{$\lambda\in \C\otimes_\Z K_0(A)$}.
\end{thm}

\begin{proof}
By the theorem of H\"ubner and Lenzing, Theorem~\ref{t:hltheorem},
there is an equivalence $\Phi:\Mod(A) \to \F^\perp$
for some finite-dimensional hereditary algebra $A$.
We denote by $D$ the duality $\Hom_\C (-,\C)$.
The endomorphism algebra of $\Phi(DA)$
is identified with the opposite algebra of $A$,
and for a coherent sheaf $\E$ on $\XX$, this turns $D\Hom(\E,\Phi(DA))$ into a left $A$-module.
We show that the functor $\ell':\Coh\XX\to \Mod(A)$ defined by
\[
\ell'(\E) = D\Hom(\E,\Phi(DA)).
\]
is left adjoint to the composition
\[
\Mod(A) \xrightarrow{\Phi} \F^\perp \xrightarrow{j_\XX} \Coh\XX.
\]
For $M\in \Mod(A)$ and $E\in \Coh\XX$ there is a map
\[
\Hom(\E,\Phi(M)) \to \Hom_A(\Hom_A(M,DA),\Hom(\E,\Phi(DA))).
\]
sending $\theta\in \Hom(\E,\Phi(M))$ to the map which sends $\alpha\in \Hom_A(M,DA)$ to $\Phi(\alpha)\theta$.
Identifying $\Hom_A(M,DA)$ with $DM$, the right hand side becomes
\[
\Hom(DM,\Hom(\E,\Phi(DA))) \cong \Hom(D\Hom(\E,\Phi(DA)),M) = \Hom(\ell'(\E),M),
\]
so we obtain a mapping
\[
t_{\E M} : \Hom(\E,\Phi(M)) \to \Hom_A(\ell'(\E),M)
\]
which is functorial in $\E$ and $M$.
Now $t_{\E M}$ is an isomorphism for $M=DA$, so it is an isomorphism
for all injective modules $M$, and then by considering the injective resolution
of a module, one sees that $t_{\E M}$ is an isomorphism for all modules $M$.
Thus $\ell'$ is left adjoint to $j_\XX\Phi$.
%Thus $\ell = \Phi \ell'$ is left adjoint to $j_\XX$.

Given $\E\in \F^\perp$, we have defined an exact sequence
\[
\mathcal{B}_\zeta(\E) : 0 \to \E(\oom) \to \mathbf{D}_\zeta(\E) \to \E \to 0.
\]
In particular, for $M\in\Mod(\Lambda)$ we have an exact sequence $\mathcal{B}_\zeta(\Phi(M))$.
Twisting, this gives an exact sequence
\[
\mathcal{B}_\zeta(\Phi(M))(-\oom) : 0 \to \Phi(M) \to \mathbf{D}_\zeta(\Phi(M))(-\oom) \to \Phi(M)(-\oom) \to 0.
\]
Applying the functor $\ell'$ to this sequence, and identifying $\ell'(\Phi(M))$ with $M$,
one obtains a right exact sequence
\[
\eta_M : M \to \ell'(\mathbf{D}_\zeta(\Phi(M))(-\oom)) \to \ell'(\Phi(M)(-\oom)) \to 0.
\]

The right hand term in $\eta_M$ is given by the assignment
\[
M \leadsto \ell'(\Phi(M)(-\oom))
\]
which defines a functor $S$ from $\Mod(\Lambda)$ to itself, and for any $M'$,
\[
\Hom(S(M),M') = \Hom(\ell'(\Phi(M)(-\oom)),M') \cong \Hom(\Phi(M)(-\oom),j_\XX \Phi(M'))
\]
\[
\cong D\Ext^1(\Phi(M'),\Phi(M)) \cong D\Ext^1(M',M)
\]
which shows that $S$ is isomorphic to the inverse Auslander-Reiten translate $\tau^-$ for $A$, which is the
functor $\tau^-M = D\Ext^1(DA,M)$ since $A$ is hereditary, and $\tau^- M \cong \tau^- A\otimes_A M$.

The middle term in $\eta_M$ is given by the assignment
\[
M \leadsto \ell'(\mathbf{D}_\zeta(\Phi(M))(-\oom))
\]
is a right exact functor from $\Mod(\Lambda)$ to itself, so it is
naturally isomorphic to $B\otimes_A M$ for some $A$-$A$-bimodule $B$.

Thus $\eta_M$ can be identified with an exact sequence
\[
\eta_M' : M \to B\otimes_A M \to \tau^-A \otimes_A M \to 0.
\]
Moreover, since this is functorial in $M$, it is obtained from a right exact sequence of $A$-$A$-bimodules
\[
A \to B \to \tau^- A \to 0.
\]
As a sequence of left $A$-modules, this is identified with the sequence $\eta_A$, and the formula for $\ell'$ shows
that it can be extended to a long exact sequence
\[
\dots \to D\Ext^1(\Phi(A)(-\oom),\Phi(DA)) \to
\]
\[
\to M \to \ell'(\mathbf{D}_\zeta(\Phi(M))(-\oom)) \to \ell'(\Phi(M)(-\oom)) \to 0.
\]
By Serre duality the first term above becomes $\Hom(\Phi(DA),\Phi(A)) \cong \Hom(DA,A)$, and the assumption that $A$ has no
nonzero projective-injective ensures that this is zero. Thus we have a sequence of bimodules
\[
0 \to A \to B \to \tau^- A \to 0
\]
which is exact on the left as well.

Now $A$ is hereditary, and the base field $\C$ is algebraically closed, so $A$ is quasi-free.
There is an exact sequence of bimodules
\[
0\to \Omega^1 A \to A\otimes_\C A \to A \to 0,
\]
and this is a projective resolution of $A$ as an $A$-$A$-bimodule, that is, as a left $A^e$-module. We write it as
\[
0\to P_1\to P_0\to A\to 0.
\]
Now the functor $M\to M^\vee = \Hom_{A^e}(M,A^e)$ gives a duality
between finitely generated projective left and right $A^e$-modules.
It gives a long exact sequence
\[
0\to A^\vee \to P_0^\vee \to P_1^\vee \to \Ext^1_{A^e}(A,A^e) \to 0
\]
This can be rewritten as
\[
0\to A^\vee \to A\otimes_\C A\to \Der(A,A\otimes_\C A) \to H^1(A,A\otimes_\C A)\to 0.
\]
Since $\Der(A,A\otimes_\C A)$ is a projective $A^e$-module, one gets a lifting
$\Der(A,A\otimes_\C A)\to B$, so a commutative diagram of $A$-$A$-bimodules
\[
\begin{CD}
& & A\otimes_\C A @>>> \Der(A,A\otimes_\C A) @>>> H^1(A,A\otimes_\C A) @>>> 0 \\
& & @VVV @VVV @| \\
0 @>>> A @>>> B @>>> \tau^- A @>>> 0
\end{CD}
\]
with exact rows, for some map $A\otimes_\C A\to A$.
Let $\alpha\in A$ be the image of $1\otimes 1$.
The fact that the lower sequence is exact on the left ensures
that this is a pushout diagram.

The category we are interested in has objects the coherent sheaves $\E\in\F^\perp$ equipped
with a $\zeta$-connection, so a section of the exact sequence $\mathcal{B}_\zeta(\E)$.
Equivalently, its objects are $A$-modules $M$ and sections of the exact sequence
\[
\mathcal{B}_\zeta(\Phi(M)) : 0 \to \Phi(M)(\oom) \to \mathbf{D}_\zeta(\Phi(M)) \to \Phi(M) \to 0.
\]
Such sections are in 1-1 correspondence with retractions, and so also in 1-1 correspondence with
retractions of the twisted sequence
\[
\mathcal{B}_\zeta(\Phi(M))(-\oom) : 0 \to \Phi(M) \to \mathbf{D}_\zeta(\Phi(M))(-\oom) \to \Phi(M)(-\oom) \to 0.
\]
Applying the functor $\ell'$ to this sequence, and identifying $\ell'(\Phi(M))$ with $M$,
one obtains a right exact sequence
\[
\eta_M : M \to \ell'(\mathbf{D}_\zeta(\Phi(M))(-\oom)) \to \ell'(\Phi(M)(-\oom)) \to 0.
\]
Now retractions $r$ of $\mathcal{B}_\zeta(\Phi(M))(-\oom)$ are in 1-1 correspondence
with retractions of $\eta_M$, that is, mappings
\[
r' : \ell'(\mathbf{D}_\zeta(\Phi(M))(-\oom)) \to M
\]
such that the composition
\[
M \to \ell'(\mathbf{D}_\zeta(\Phi(M))(-\oom)) \to M
\]
is the identity map on $M$.
Namely, given $r$ one obtains $r'$ by applying $\ell'$ to $r$, and given $r'$, one
obtains $r$ by composing $j_\XX \Phi (r')$ with the adjunction map
\[
\mathbf{D}_\zeta(\Phi(M))(-\oom) \to j_\XX \Phi \ell'(\mathbf{D}_\zeta(\Phi(M))(-\oom)).
\]
Now retractions $r'$ are in 1-1 correspondence with retractions of the sequence
\[
\eta_M' : M \to B\otimes_A M \to \tau^-A \otimes_A M \to 0.
\]
Now we have a pushout diagram
\[
\begin{CD}
A\otimes_\C A @>>> \Der(A,A\otimes_\C A) \\
@VVV @VVV \\
A @>>> B
\end{CD}
\]
and it remains a pushout after tensoring with $M$. Thus a retraction of $M\to B\otimes_A M$ is the same thing
as an $A$-module map $\Der(A,A\otimes_\C A)\otimes_A M \to M$ which sends $\Delta\otimes m$ to $\alpha m$,
where $\Delta$ is the derivation with $\Delta(a)=a\otimes 1 - 1\otimes a$.
This is the same as the data for a module for the preprojective algebra $\Pi^\alpha(A)$.
Finally, since $A$ is a finite-dimensional hereditary algebra, the
mapping $\C\otimes_\Z K_0(A)\to A/[A,A]$ used in \cite{CBconze} is an isomorphism
(as one easily sees by using Morita equivalence to pass to the corresponding basic
algebra, and using that this is isomorphic to the path algebra of a quiver without oriented cycles).
Thus $\Pi^\alpha(A) \cong \Pi^\lambda(A)$ for some $\lambda\in \C\otimes_\Z K_0(A)$.
\end{proof}

\end{document}